\documentclass[a4paper,10pt]{article}

\usepackage{amsmath}
\usepackage{ntheorem}
\usepackage{eucal}
\usepackage{amscd}
\usepackage{ifthen}
\usepackage{amssymb}
\usepackage[english,french]{babel}
\usepackage[T1]{fontenc}
\usepackage[latin1]{inputenc}
\usepackage[all]{xy}

\newenvironment{demo}[1][]{\ifthenelse{\equal{#1}{}}{\noindent \textbf{D\'emonstration :\\ \indent}}{\noindent \textbf{D\'emonstration #1 :\\ \indent}}}{$\square$\\}

\theoremstyle{break}
\newtheorem{theo}{Théorème}[section]
\newtheorem{cor}{Corollaire}[section]
\newtheorem{conj}{Conjecture}
\newtheorem{lem}{Lemme}[section]
\newtheorem{prop}{Proposition}[section]
\newtheorem{defi}{Définition}[section]
\theorembodyfont{\slshape}
\newtheorem{rem}{Remarque}[section]
\theorembodyfont{\upshape}

\newcommand{\CC}{\mathbb C}
\newcommand{\QQ}{\mathbb Q}

\newcommand{\DD}{\mathbb D}
\newcommand{\PP}{\mathbb P}
\newcommand{\XX}{\mathcal{X}}
\newcommand{\kod}[1]{\kappa(#1)}
\newcommand{\gd}[1]{\gamma d(#1)}
\newcommand{\To}{\longrightarrow}
\newcommand{\abs}[1]{\left\vert#1\right\vert}

\newcommand{\Forall}[2]{\forall \, #1 \in #2, \:}

\newcommand{\cycl}[1]{\mathcal{C}(#1)}
\newcommand{\hg}[2]{\pi_{#1}(#2)}
\newcommand{\dimm}[1]{\mathrm{dim}(#1)}
\newcommand{\alb}[1]{\mathrm{Alb}(#1)}

\newcommand{\vol}[1]{\mathrm{Vol}(#1)}
\newcommand{\wtx}{\widetilde{X}}

\title{Invariance de la $\Gamma$-dimension pour certaines familles kählériennes de dimension 3}%[Invariance de la $\Gamma$-dimension en dimension 3]
\author{Benoît Claudon}
\date{}

\begin{document}

\maketitle

\begin{abstract}
Dans cet article, on étudie les propriétés d'invariance par déformation de la $\Gamma$-dimension d'une variété kählérienne compacte $X$, cet invariant birationnel contrôlant la dimension des sous-variétés compactes maximales contenues dans le revêtement universel de $X$. Pour les surfaces kählériennes, cette propriété d'invariance résulte immédiatement d'un théorème de Y.-T. Siu. En utilisant un théorème de structure de F. Campana et Q. Zhang, on montre également l'invariance par déformation de la $\Gamma$-dimension pour certaines familles de variétés kählériennes compactes de dimension 3.
\end{abstract}
\vspace*{0.5cm}
\begin{center}
{\LARGE Invariance of $\Gamma$-dimension for some Kähler\\
\vspace*{0.4cm}3-folds families}
\end{center}
\selectlanguage{english}
\begin{abstract}
In this article, we study some properties of deformation invariance of the $\Gamma$-dimension (defined for $X$ a compact kähler manifold). This birational invariant is defined as the codimension of the maximal compact subvarieties in the universal cover of $X$. In the surface case, the deformation invariance is a straightforward consequence of a theorem of Y.-T. Siu. Using some results from F. Campana et Q. Zhang, we settle this invariance for certain type of Kähler families of dimension 3.
\end{abstract}
\selectlanguage{french}

\section*{Introduction}

La théorie des fonctions, ainsi que la structure du revêtement universel d'une variété projective (ou plus généralement kählérienne compacte), sont conjecturalement décrites par l'énoncé suivant \cite{Sh74} :
\begin{conj}[I.R. Shafarevich, 1974]\label{shaf1}
Soit $X$ une variété projective de revêtement universel $\wtx$. La variété complexe $\wtx$ est alors holomorphiquement convexe ; en particulier, $\wtx$ admet une application propre sur un espace de Stein (appelée réduction de Remmert de $\wtx$).
\end{conj}

A défaut de pouvoir construire \emph{a priori} des fonctions holomorphes non constantes sur $\wtx$, on peut construire une réduction de Remmert pour $\wtx$ en un sens biméromorphe :
\begin{theo}[\cite{Ca94}]\label{existence-gred}
Si $X$ est une variété kählérienne compacte, $\wtx$ admet une unique application méromorphe, presque-holomorphe\footnote{cela signifie que le lieu d'indétermination ne s'envoie pas surjectivement sur la base de la fibration.}, propre et à fibres connexes
$$\gamma_{\wtx}:\wtx\dashrightarrow\Gamma(\wtx)$$
vérifiant : si $\widetilde{Z}$ est une sous-variété irréductible compacte qui intersecte une fibre générale de $\gamma_{\wtx}$, alors $\widetilde{Z}$ est entièrement contenue dans celle-ci.\\
De façon équivalente, il existe sur $X$ une unique fibration presque-holomorphe
$$\gamma_X:X\dashrightarrow \Gamma(X)$$
vérifiant : si $Z$ est une sous-variété irréductible passant par un point général $x\in X$, de normalisée $\hat{Z}\To Z$ et satisfaisant à la condition
$$\abs{\mathrm{Im}(\hg{1}{\hat{Z}}\To\hg{1}{X}}<+\infty$$
$Z$ est alors contenue dans la fibre en $x$
$$Z\subset \gamma_X^{-1}(\gamma_X(x)).$$
\end{theo}
\begin{rem}
La conjecture \ref{shaf1} est alors équivalente au fait que $\gamma_{\wtx}$ soit holomorphe et que $\Gamma(\wtx)$ soit de Stein.
\end{rem}
\begin{rem}\label{remarque kollar}
Indépendament, J. Koll\'ar a également proposé une construction de $\gamma_X$ dans le cas algébrique. Dans \cite{K93}, elle est notée $Sh_X$ et prend le nom d'application de Shafarevich.
\end{rem}
Dans la suite, on notera $\hg{1}{Z}_X$ l'image du morphisme $\hg{1}{Z}\To \hg{1}{X}$ pour une sous-variété $Z$ de $X$ et $\hg{1}{\hat{Z}}_X$ l'image du morphisme induit par la normalisation de $Z$.\\

\noindent Le théorème \ref{existence-gred} permet alors de définir l'invariant suivant.
\begin{defi}\label{defi-gred}
On appelle $\Gamma$-dimension de la variété $X$ l'entier défini par :
$$\gd{X}=\dimm{\Gamma(X)}$$
C'est un invariant birationnel qui ne dépend que du revêtement universel de $X$ (il est en particulier invariant sous les revêtements étales finis de $X$).
\end{defi}

La question de l'invariance par déformation (qui se pose naturellement) est soulevée par J. Koll\'ar dans \cite[18.4, p. 184]{K95} et étudiée dans \cite{OKR}. Formulons la explicitement :
\begin{conj}\label{conjecture-déformation}
La $\Gamma$-dimension est invariante dans une famille de variétés kählériennes\footnote{voir la section suivante pour les différentes notions de familles utilisées ici.}.
\end{conj}
Dans l'article \cite{OKR}, les auteurs s'intéressent notamment à l'invariance par déformation pour les variétés\footnote{ces variétés sont dites de type $\pi_1$-général ; voir la section \ref{section pi1-general}.} $X$ vérifiant :
$$\gd{X}=\dimm{X}.$$

Les résultats principaux développés dans les sections suivantes traitent le problème de l'invariance par déformation de la $\Gamma$-dimension pour certaines familles de dimension 3. Ainsi, on démontrera dans cet article les théorèmes suivants :
\begin{theo}\label{invariance non type general}
Soit $\pi:\XX\To B$ une famille de variétés kählériennes compactes de dimension 3. Si la famille n'est pas de type général (\emph{i.e.} si les fibres ne le sont pas), la $\Gamma$-dimension est constante dans la déformation : l'application
$$B\ni b\mapsto \gd{\XX_t}$$
est localement constante sur $B$.
\end{theo}
De plus, sous réserve d'un résultat conjectural, la condition $\gd{X}=3$ est une condition ouverte dans une déformation :
\begin{theo}[énoncé conditionnel]\label{ouverture gd=3}
Soit $\pi:\XX\To \DD$ une famille projective de variétés de dimension 3. Si la fibre centrale $\XX_0$ vérifie $\gd{\XX_0}=3$, les fibres voisines de $\XX_0$ sont alors également de type $\pi_1$-général :
$$\gd{\XX_t}=3\textrm{ pour tout $t$ dans un voisinage de }0.$$
\end{theo}
La conjecture dont dépend le théorème précédent peut se formuler de la façon suivante :
\begin{conj}\label{conjecture surface}
Si $S$ est une surface rationnelle et $\Delta$ un $\QQ$-diviseur (dont le support est à croisements normaux) de la forme :
$$\Delta=\sum_{j\in J}(1-\frac{1}{m_j})\Delta_j \textrm{ (avec $m_j\geq2$ des entiers)}$$
tels que $\kod{S/\Delta}\leq 0$, le groupe fondamental orbifolde\footnote{voir la section \ref{section gd=2}.} $\hg{1}{S/\Delta}$ est presque abélien (\emph{i.e.} admet un sous-groupe abélien d'indice fini).
\end{conj}

Ce papier est structuré de la façon suivante : la première section est dévolue à quelques rappels sur les familles de variétés kählériennes, l'existence d'une $\Gamma$-réduction générique dans une telle famille ainsi qu'à l'exposition du cas des surfaces kählériennes. La deuxième partie aura pour objectif de montrer le théorème \ref{invariance non type general} ; pour cela, nous nous appuyerons sur un théorème de structure dû à F. Campana et Q. Zhang. Enfin, nous étudierons plus précisément les variétés de dimension 3 vérifiant $\gd{X}=2$ pour en déduire le théorème \ref{ouverture gd=3}.

\section{Propriétés relatives de la $\Gamma$-réduction}

\subsection{Structure locale des familles de variétés kählériennes}

On rappelle ici quelques notions de base sur les familles de variétés et sur la structure locale de celles-ci lorsqu'elles sont suposées kählériennes. Dans un souci de complétude, nous redonnons également la démonstration d'un énoncé de compacité relative dans l'espace des cycles (pour les notions concernant l'espace des cycles et sa topologie, on renverra systématiquement à \cite{Ba75} et \cite{Li78}).
\begin{defi}\label{defi-famille}
Une famille de variétés complexes est un triplet $(\XX,B,\pi)$ dans lequel $\XX$ et $B$ sont des variétés complexes et
$$\pi:\XX\To B$$
une submersion holomorphe propre dont les fibres sont des variétés connexes (et compactes par propreté de $\pi$). La famille $(\XX,B,\pi)$ est dite :
\begin{enumerate}
\item {\bf kählérienne} si toutes les fibres de la famille le sont.
\item {\bf projective} s'il existe un fibré en droites $A$ sur $\XX$ qui est $\pi$-ample.
\end{enumerate}
\end{defi}
\begin{rem}
Si $\pi:\XX\To B$ est une famille projective au-dessus d'une base Stein, on peut toujours se ramener (et c'est ce que l'on fera dans la suite) au cas où le fibré $A$ est ample sur $\XX$.
\end{rem}

La structure locale (sur la base) d'une famille de variétés kählériennes est en fait assez simple : on dispose d'une métrique "relativement kählérienne".
\begin{prop}\label{metrique-locale}
Soit $\pi:\XX\To B$ une famille de variétés kählériennes compactes et $b$ un point de $B$. Il existe alors un voisinage (que l'on peut supposer contractile) $U$ de $b$ dans $B$ et $\omega$ une métrique hermitienne sur $X=\XX_U=\pi^{-1}(U)$ qui se restreint en une métrique kählérienne sur chaque fibre de la famille $\pi:X\To U$ ($i.e.$ $d(\omega_{\vert X_u})=0$ pour tout $u\in U$).
\end{prop}
\begin{demo}
On peut facilement s'en convaincre en examinant la démonstration du théorème de Kodaira sur les petites déformations des variétés kählériennes, par exemple celle proposée dans \cite[th. 9.23, p. 222]{V02}. On y construit la métrique pour montrer que les fibres voisines de $\XX_b$ sont kählériennes.
\end{demo}

\noindent Dans cette situation, les propriétés de compacité relative des cycles sont préservées comme le montre la proposition suivante (en comparaison de la situation absolue où $U$ est réduit à un point \cite{Li78}).
\begin{prop}\label{propreté-famille}
Si $\pi:X\To U$ est une famille de variétés kählériennes compactes comme ci-dessus ($U$ est contractile et $X$ est muni d'une métrique $\omega$ qui se restreint en une métrique kählérienne sur chaque fibre), les composantes irréductibles de $\cycl{X/U}$ sont alors $U$-propres pour la projection
$$\pi_*:\cycl{X/U}\To U$$
\end{prop}
Rappelons que sans l'hypothèse contractile, la propriété de compacité relative ci-dessus est mise en défaut comme le montre l'exemple de D. Lieberman ($ibid$). Nous redonnons la démonstration de cette proposition faute de références dans la littérature existante.\\

\begin{demo}
Comme $U$ est contractile, le lemme d'Ehresmann (voir \cite[prop. 9.3, p. 209]{V02}) montre alors que la fibration est topologiquement triviale :
$$ \xymatrix{
X \ar[rr]^{\varphi} \ar[rd]_{\pi} && X_0\times U \ar[ld]^{pr_2}\\
 & U}$$
où $\varphi$ est un difféomorphisme et $0\in U$ un point base quelconque. En particulier, les algèbres de cohomologie $H^\bullet(X_u)$ (pour $u\in U$) sont canoniquement isomorphes à $H^\bullet(X)$. Si $\omega_u$ désigne la restriction de $\omega$ à la fibre $X_u$, $\omega_u$ définit une classe dans $H^2(X_u)$ (car elle est $d$-fermée) et elle définit donc également une classe $[\omega_u]\in H^2(X)$.

Soit maintenant $(C_t)_{t\in T}$ une famille analytique de $n$-cycles (relatifs) paramétrée par un espace irréductible $T$ et supposons que
$$\Forall{t}{T}u(t)=\pi_*(C_t)\in K$$
avec $K$ un compact de $U$. Le volume est alors donné par :
$$\textrm{Vol}_\omega(C_t)=\int_{C_t}\omega^n= \int_{C_t}(\omega_{u(t)})^n=[C_t]\wedge[\omega_{u(t)}]^n$$
où le produit d'intersection est calculé (par exemple) dans $H^\bullet(X)$. Mais lorsque $u$ décrit $K$, $[\omega_u]$ décrit une partie bornée de $H^2(X)$ (par compacité de $K$) et il est bien connu que dans, une famille irréductible, la classe d'homologie est fixée : on en déduit que le volume de la famille $\textrm{Vol}_\omega(C_t)$ est borné. D'autre part, les supports des cycles $C_t$ sont contenus dans le compact $\pi^{-1}(K)$. On peut alors appliquer \cite{Li78} et conclure que la famille $(C_t)_{t\in T}$ est contenue dans un compact de $\cycl{X}$ ; comme $\cycl{X/U}$ est fermé dans $\cycl{X}$, $(C_t)_{t\in T}$ est contenue dans un compact de $\cycl{X/U}$ et cela montre la propreté de $\pi_*$ en restriction aux composantes irréductibles.
\end{demo}

\subsection{$\Gamma$-réduction relative}

Rappelons tout d'abord que la construction de la $\Gamma$-réduction peut se faire en famille mais, \emph{a priori}, seulement génériquement :
\begin{prop}\label{gred-famille-générique}
Soit $\pi:\XX\To B$ une famille de variétés kählériennes compactes. Il existe alors $\Gamma(\XX/B)$ un espace complexe normal et un diagramme :
$$\xymatrix{\XX\ar@{-->}[rr]^{\gamma_{\pi}}\ar[rd]_{\pi} &  & \Gamma(\XX/B)\ar[ld]^{\phi}\\
& B &}$$
dans lequel $\phi$ est surjectif et telle que pour $b\in B$ général, la restriction de $\gamma_{\pi}$ à $\XX_b$ :
$$\gamma_{\pi\vert \XX_b}:\XX_b\dashrightarrow \Gamma(\XX/B)$$
soit la $\Gamma$-réduction de $\XX_b$.
\end{prop}
\begin{demo}
On applique les théorèmes 2.2 et 2.7 de \cite{Ca04app}. On pourra également comparer à \cite[rk. 2.7]{OKR}
\end{demo}

On peut maintenant, grâce à la proposition \ref{propreté-famille} ci-dessus, montrer une propriété de semi-continuité inférieure de la $\Gamma$-dimension au cours d'une déformation.
\begin{theo}\label{gd-générique=max}
Si $\pi:\XX\To B$ est une famille de variétés kählériennes compactes et si $d$ désigne la valeur générale de la $\Gamma$-dimension de la famille $\pi$, alors
$$d=\underset{b\in B}{\mathrm{max}}\big(\gd{\XX_b}\big)$$
\end{theo}
\begin{demo}
Notons $n$ la dimension relative de $\XX$ sur $B$ et fixons $b$ un point de $B$ et $U$ un voisinage contractile (par exemple un polydisque) de $b$ comme dans la proposition \ref{metrique-locale}. On notera encore $\XX$ la restriction de la famille $\pi$ au dessus de $U$. Comme l'espace $\mathcal{C}_{n-d}(\XX/U)$ n'a qu'un nombre au plus dénombrable de composantes irréductibles, il en existe une que l'on notera $\mathfrak{C}$ qui contient les fibres de $\gamma_{\XX_u}$ pour $u$ général dans $U$. Or, d'après la proposition \ref{propreté-famille}, l'application
$$\pi_*:\mathfrak{C}\To U$$
est propre donc en particulier fermée. D'autre part, $\pi_*(\mathfrak{C})$ contient le sous-ensemble dense des points généraux ci-dessus ; on en déduit donc que $\pi_*$ est surjective. Au dessus du point $b$, $\mathfrak{C}_b=\pi_*^{-1}(b)$ est donc un sous-ensemble analytique de $\mathcal{C}_{n-d}(\XX_b)$ dont une composante irréductible au moins doit être couvrante (pour $\XX_b$). Il nous reste maintenant à vérifier l'affirmation suivante : les (composantes irréductibles des) cycles $(Z_t)_{t\in T}$ de dimension $n-d$ paramétrés par cette composante irréductible vérifient donc
$$\textrm{Im}(\hg{1}{\hat{Z_t}}\To\hg{1}{\XX_b})<+\infty$$
En effet, si on admet le fait ci-dessus, on en déduit immédiatement (par définition de l'invariant $\gamma d$) :
$$\Forall{b}{B}\gd{\XX_b}\leq d$$
On va en fait montrer que les composantes irréductibles des images réciproques (dans $\widetilde{\XX_b}$) des cycles $(Z_t)_{t\in T}$ sont compactes (ce qui est équivalent à l'assertion ci-dessus). Si on se place au dessus d'un voisinage compact $K$ de $b$, le volume des cycles de $\mathfrak{C}_K=\pi_*^{-1}(K)$ par rapport à la métrique $\omega$ de la proposition \ref{metrique-locale} est uniformément majoré :
$$\Forall{Z}{\mathfrak{C}_K}\mathrm{Vol}_{\omega}(Z)\leq M$$
par une constante $M>0$. Pour les cycles $Z$ qui sont les fibres des $\Gamma$-réduction de $\XX_u$ ($u\in K$ général et $u\neq b$), on a donc également :
$$\mathrm{Vol}_{\tilde{\omega}}(Z^j)\leq M'$$
où $\tilde{\omega}$ est l'image réciproque de $\omega$ sur $\widetilde{\XX}$ (le revêtement universel de $\XX$) et $\tilde{Z}^j$ une composante irréductible de l'image réciproque de $Z$ dans $\widetilde{\XX}$. La constante $M'$ vérifie en fait $M'=NM$ avec
$$N=\textrm{Card}\left(\hg{1}{Z}_{\XX}\right)$$
(cet entier ne dépend pas de $Z$ générique dans la famille $\mathfrak{C}_K$). En passant à la limite, on en déduit que les sous-variétés $\tilde{Z_t}^j$ ($t\in T$) de $\widetilde{\XX_b}$ ont un volume fini ; un argument standard de géométrie bornée (voir \cite[p. 288-289]{G91}) montre alors que les $\tilde{Z_t}^j$ sont compactes.
\end{demo}

\subsection{Invariance dans le cas des surfaces}

Rappelons pour finir que la $\Gamma$-dimension est effectivement invariante dans les familles de surfaces kählériennes mais ce pour des raisons tenant au fait que la dimension est petite. En effet, d'après un théorème de Y.-T. Siu, la propriété $\gd{X}=1$ ne dépend que du groupe fondamental de $X$ :
\begin{theo}[\cite{S87}]\label{theoreme-siu}
Si $X$ est une variété kählérienne compacte, alors $\gd{X}=1$ si et seulement si $\hg{1}{X}$ est commensurable\footnote{on rappelle que deux groupes $G$ et $H$ sont dits commensurables si ils possèdent des sous-groupes d'indices finis $G_1\subset G$ et $H_1\subset H$ isomorphes.} au groupe fondamental d'une courbe $C$ de genre $g\geq 1$.
\end{theo}
Les conditions $\gd{X}=0$ et $\gd{X}=1$ ne dépendent donc que du groupe fondamental de $X$.
\begin{cor}\label{gdim-surface}
Si $X$ est une surface kählérienne compacte, la $\Gamma$-dimension de $X$ ne dépend que du groupe fondamental (c'est en particulier un invariant topologique de $X$).
\end{cor}
Notons que ceci n'est bien entendu plus valable en dimension supérieure. Considérons en effet $A$ une variété projective lisse de dimension $n\geq 4$ telle que $\gd{A}=n$, $X$ une section hyperplane de $A$ et $Z$ une section hyperplane de $X$. Si $Y$ est obtenue comme section hyperplane de $\PP^2\times Z$, on vérifie facilement que $X$ et $Y$ sont deux variétés projectives de dimension $n-1$ vérifiant :
$$\gd{X}=n-1\textrm{ et }\gd{Y}=n-2.$$
Or, des application répétées du théorème de Lefschetz montrent que
$$\hg{1}{X}\simeq\hg{1}{Y}\simeq\hg{1}{A}.$$
Ainsi, $X$ et $Y$ sont deux variétés projectives (de dimension au moins 3) ayant des groupes fondamentaux isomorphes mais de $\Gamma$-dimension différentes.
\begin{cor}\label{invariance-gdim-surface}
Si $(\XX,B,\pi)$ est une déformation de surfaces kählériennes au dessus d'une base $B$ connexe, la $\Gamma$-dimension est invariante dans la famille. En d'autres termes, l'application
$$B\ni b\mapsto \gd{\XX_b}$$
est constante sur $B$.
\end{cor}
\begin{demo}[du corollaire \ref{invariance-gdim-surface}]
D'après le corollaire \ref{gdim-surface}, la $\Gamma$-dimension ne dépend que du groupe fondamental pour une surface kählérienne compacte. Or, dans une famille $(\XX,B,\pi)$, les fibres sont deux à deux homéomorphes d'où l'invariance par déformation.
\end{demo}

Dans les cas des familles de dimension 3, nous sommes donc face à la dichotomie suivante : $\gd{X}=2$ ou $3$ ; dichotomie que nous allons étudier dans les sections suivantes.

\section{Variétés de type $\pi_1$-général}\label{section pi1-general}

\subsection{Un résultat de structure}

Rappelons que les variétés $X$ vérifiant $\gd{X}=\dimm{X}$ sont dites de type $\pi_1$-général (voir \cite{CZ04} et \cite[th. 4.1]{Ca95} pour la justification de cette terminologie). La structure de ces variétés est conjecturalement décrite par l'énoncé suivant.
\begin{conj}[\cite{K93}, \cite{CZ04}]\label{conjecture type pig}
Une variété kählérienne compacte $X$ de type $\pi_1$-général vérifie nécessairement $\kod{X}\geq 0$ et, à revêtement étale fini et transformation birationnelle près, $X$ admet une submersion en tores $f:X\To Y$ sur une variété $Y$ de type général et de type $\pi_1$-général.
\end{conj}
Dans le cas projectif, cette conjecture est une conséquence de la conjecture d'Abondance et des résultats de Bogomolov-Yau sur les variétés à fibré canonique trivial. Nous rappelons ici le résultat principal de \cite{CZ04} qui établit la conjecture \ref{conjecture type pig} en dimension 3.
\begin{theo}\label{CZ}
Soit $X$ une variété kählérienne compacte de dimension 3 avec $\gd{X}=3$. On a alors $\kod{X}\geq 0$ et de plus, à revêtement étale fini près, $X$ est biméromorphe
\begin{enumerate}
\item soit à $\mathrm{Alb}(X)\times_{\mathrm{Alb}(J(X))}J(X)$
\item soit à une submersion en surfaces abéliennes sur une courbe de genre au moins 2.
\end{enumerate}
En particulier, toujours à revêtement étale fini près, la fibration d'Iitaka-Moishezon $J_X:X\To J(X)$ de $X$ est (biméromorphe à) une submersion en tores et $J(X)$ est de type général et de type $\pi_1$-général.
\end{theo}
Les propositions suivantes (immédiates avec le théorème ci-dessus) isolent les cas $\kod{X}=1$ et 2.
\begin{prop}\label{gd3/kappa1}
Soit $X$ kähler compacte de dimension 3 vérifiant $\gd{X}=3$ et $\kod{X}=1$. A revêtement étale fini près, $X$ admet une submersion en tores sur une courbe
$$J:X\To C$$
vérifiant $g(C)\geq 2$ et, si $T$ désigne la fibre de $J$, le morphisme $\hg{1}{T}\To \hg{1}{X}$ est injectif.
\end{prop}
\begin{prop}\label{gd3/kappa2}
Soit $X$ kähler compacte de dimension 3 vérifiant $\gd{X}=3$ et $\kod{X}=2$. A revêtement étale fini près, l'application d'Albanese de $X$ est
\begin{enumerate}
\item soit isotriviale sur une courbe elliptique avec pour fibre une surface de type $\pi_1$-général.
\item soit une fibration sur une surface de type $\pi_1$-général dont la fibre générale $F$ vérifie $\hg{1}{F}_X$ est un groupe infini (avec la notation introduite après la remarque \ref{remarque kollar}).
\item soit génériquement finie sur son image.
\end{enumerate}
\end{prop}

Pour finir ce paragraphe, rappelons les faits suivants (qui nous seront utiles ultérieurement) concernant l'application d'Albanese. Si on note
$$\alpha_X:X\To \alb{X}$$
l'application d'Albanese de la variété $X$ et $\alpha d(X)=\dimm{\alpha_X(X)}$ sa dimension d'Albanese, on dispose alors des résultats suivants.
\begin{prop}[prop. 1.4, p. 267 \cite{Cat91}]\label{inv-albdim}
La dimension d'Albanese est un invariant topologique (et même un invariant d'homotopie).
\end{prop}
\begin{prop}\label{gdim-albdim}
Soit $X$ une variété kählérienne compacte. On a alors
$$\gd{X}\geq \alpha d(X).$$
\end{prop}
\begin{demo}
L'image de $\alpha_X$ est une sous-variété du tore $\alb{X}$ et elle est donc de type $\pi_1$-général. Comme la $\Gamma$-dimension est croissante sous les applications surjectives, on a bien
$$\gd{X}\geq \alpha d(X).$$
\end{demo}

\subsection{$\Gamma$-réduction et fibrations}

Dans ce paragraphe, nous établissons une sorte de réciproque à la conjecture \ref{conjecture type pig}.
\begin{prop}\label{submersion tore}
Si une variété kählérienne compacte $X$ admet une submersion en tores $f:X\To Y$, on a alors :
$$\gd{X}\geq\dimm{f}+\gd{Y}.$$
En particulier, $X$ est de type $\pi_1$-général si $Y$ l'est.
\end{prop}
Le comportement de la $\Gamma$-dimension dans une fibration $f:X\To Y$ entre variétés kählériennes compactes est en effet assez difficile à décrire en toute généralité. On sait en effet que $\gd{X}\geq \gd{Y}$ mais, même si la fibre générale $F$ de $f$ vérifie $\gd{F}>0$, l'inégalité ci-dessus n'est pas nécessairement stricte. Par exemple, si $X$ est une section hyperplane de $\PP^2\times\PP^1$, les fibres de la projection sur le deuxième facteur
$$pr_2:X\To \PP^1$$
sont des courbes de genre $g>0$ (si $X$ est de degré au moins 3 sur $\PP^2$). Ces fibres vérifient donc $\gd{F}>0$ alors que $\gd{X}=0$ (par le théorème de Lefschetz, $X$ est simplement connexe).

Pour établir la proposition \ref{submersion tore}, il nous faut montrer que le morphisme naturel entre les groupes fondamentaux de la fibre et celui de $X$ est injectif, le lemme suivant permettant ensuite de conclure.
\begin{lem}\label{additivité}
Soit $f:X\To Y$ une fibration entre variétés kählériennes compactes dont on note $F$ la fibre générale. On suppose de plus que le morphisme $\hg{1}{F}\To \hg{1}{X}$ déduit de l'inclusion est injectif. On a alors :
$$\gd{Y}+\dimm{F}\geq\gd{X}\geq\gd{F}+\gd{Y}$$
En particulier, $X$ est de type $\pi_1$-général si $F$ et $Y$ le sont.
\end{lem}
\begin{demo}
Considérons le diagramme (fonctorialité de la $\Gamma$-réduction) :
$$\xymatrix{X \ar[r]^{f}\ar[d]_{\gamma_X} & Y\ar[d]^{\gamma_Y}\\
\Gamma(X)\ar[r]^{\gamma_f} & \Gamma(Y)
}$$
L'image $\gamma_X(F)$ de $F$ par $\gamma_X$ est contenue dans une fibre de
$$\gamma_f:\Gamma(X)\To \Gamma(Y)$$
d'où :
\begin{equation}\label{ineg-gred}
\dimm{\gamma_X(F)}\leq \dimm{\Gamma(X)}-\dimm{\Gamma(Y)}=\gd{X}-\gd{Y}
\end{equation}
Notons $\Gamma$ la fibre générale de la $\Gamma$-réduction  de $X$ et considérons l'application $F\To \gamma_X(F)$. La fibre de cette application est $\Gamma\cap F$ et comme $\hg{1}{F}$ s'injecte dans $\hg{1}{X}$, on constate que l'image de $\hg{1}{\Gamma\cap F}$ dans $\hg{1}{F}$ doit donc être finie. Donc $\Gamma\cap F$ est contenue dans la fibre de la $\Gamma$-réduction  de $F$ et l'application $F\To \gamma_X(F)$ se factorise sous la forme :
$$\xymatrix{F \ar[rr] \ar[rd] && \gamma_X(F)\ar[ld]\\
 & \Gamma(F)
}$$
et on en déduit que $\dimm{\gamma_X(F)}\geq\gd{F}$. Ceci combiné avec l'inégalité (\ref{ineg-gred}) fournit l'additivité annoncée :
$$\gd{X}\geq\gd{F}+\gd{Y}.$$
\end{demo}

\noindent Lorsque la fibre $F$ est une courbe, la démonstration ci-dessus fournit le corollaire immédiat suivant.
\begin{lem}\label{additivité-courbe}
Si $f:X\To Y$ est une fibration en courbes et si l'image du morphisme $\hg{1}{F}\To \hg{1}{X}$ est infinie, on a alors :
$$\gd{X}\geq 1+\gd{Y}$$
\end{lem}

Comme annoncé ci-dessus, il nous reste à montrer l'injectivité du morphisme $\hg{1}{F}\To \hg{1}{X}$ dans la situation considérée.
\begin{lem}\label{injectivité}
Si $f:X\To Y$ est une submersion en tores avec $X$ kählérienne compacte (de fibre $F$), le morphisme induit
$$i_*:\hg{1}{F}\To \hg{1}{X}$$
est injectif.
\end{lem}
\begin{rem}
L'exemple des surfaces de Hopf (au moins celles admettant une fibration elliptique sur $\PP^1$) montre que l'hypothèse "$X$ kählérienne" est essentielle.
\end{rem}
\begin{demo}
Si $\pi:\widetilde{Y}\To Y$ est le revêtement universel de $Y$, le produit fibré :
$$\xymatrix{\wtx_f=X\times_Y \widetilde{Y}\ar[d]_{\pi_f}\ar[r] & \widetilde{Y}\ar[d]^{\pi}\\
X\ar[r]^{f} & Y}$$
définit une submersion holomorphe propre $g:\wtx_f\To \widetilde{Y}$. De plus, l'application $\pi_f:\wtx_f\To X$ est un revêtement étale. La suite exacte d'homotopie associée à $g$ s'écrit :
\begin{equation}\label{homotopie-tore}
0\To K\To \hg{1}{F}\stackrel{i_*}{\To} \hg{1}{\wtx_f}\To 0
\end{equation}
où $K$ désigne simplement le noyau de $\hg{1}{F}\To \hg{1}{\wtx_f}$. La suite (\ref{homotopie-tore}) est donc une suite exacte de groupes abéliens qui se dualise en :
\begin{equation}\label{homotopie-tore-dual}
0\To H^1(\wtx_f,\CC)\stackrel{i^*}{\To } H^1(F,\CC)\To \textrm{Hom}(K,\CC).
\end{equation}
Or, comme $\widetilde{Y}$ est simplement connexe, la monodromie de la fibration $g$ est triviale et, comme $\wtx_f$ est kählérienne, on peut appliquer le théorème de dégénérescence de la suite spectrale de Blanchard \cite{Bl56} (voir également \cite[p. 376-382]{V02} pour une exposition remarquablement limpide de ce résultat) : le morphisme de restriction
$$H^1(\wtx_f,\CC)\stackrel{i^*}{\To } H^1(F,\CC)$$
est surjectif et la flèche
$$ H^1(F,\CC)\To \textrm{Hom}(K,\CC)$$
est donc identiquement nulle. Comme $K$ est un sous-groupe d'un groupe abélien sans torsion, ceci est équivalent au fait que $K$ soit réduit à zéro. On a donc un isomorphisme de groupes :
$$i_*:\hg{1}{F}\stackrel{\sim}{\To} \hg{1}{\wtx_f}.$$
Pour conclure, rappelons que le revêtement $\pi_f:\wtx_f\To X$ est étale ; le morphisme induit au niveau du groupe fondamental
$$\hg{1}{F}\stackrel{\sim}{\To} \hg{1}{\wtx_f}\stackrel{\pi_{f*}}{\To}\hg{1}{X}$$
est donc injectif.
\end{demo}
La proposition \ref{submersion tore} s'obtient alors en combinant les lemmes \ref{additivité} et \ref{injectivité}.

\subsection{Invariance pour les familles non de type général}

Pour démontrer le théorème \ref{invariance non type general}, nous allons utiliser les résultats de structure des paragraphes précédents. 

\begin{demo}[du théorème \ref{invariance non type general}]
Comme il faut montrer que la $\Gamma$-dimension est localement constante dans une déformation, on se donne
$$\pi:\XX\To \DD$$
une famille de variétés kählériennes compactes de dimension 3 au dessus du disque unité de $\CC$. Nous noterons $X=\XX_0$ la fibre centrale de la famille. Nous allons raisonner au cas par cas selon les valeurs de $\gd{X}$ et de $\kod{X}$ et montrer que, pour $t\in\DD$ proche de 0, on a $\gd{\XX_t}=\gd{X}$. Rappelons que seuls les cas $\gd{X}=2$ et $\gd{X}=3$ sont à traiter (et ce d'après le théorème \ref{theoreme-siu}).
\begin{rem}
La plupart des énoncés de structure ci-dessus sont à revêtement étale fini près. Un tel revêtement correspond à un sous-groupe d'indice fini de
$$\hg{1}{X}\simeq\hg{1}{\XX}$$
et induit donc un revêtement étale fini de $\XX$. Quitte à remplacer $\XX$ par ce revêtement, nous supposerons toujours que la variété $X$ possède les propriétés mentionnées dans les propositions \ref{gd3/kappa1} et \ref{gd3/kappa2}.
\end{rem}

\noindent\textbf{cas 1 : \mathversion{bold}$\gd{X}=3$ et $\kod{X}=0$}\mathversion{normal}\\
La variété $X$ est donc biméromorphe à un tore d'après \cite{CZ04} ; comme toute déformation d'un tore est un tore, on a bien :
$$\Forall{t}{\DD}\gd{\XX_t}=3=\gd{X}.$$
\vspace*{0.3cm}

\noindent\textbf{cas 2 : \mathversion{bold}$\gd{X}=3$ et $\kod{X}=1$}\mathversion{normal}\\
La proposition \ref{gd3/kappa1} montre que $X$ admet une submersion en tores $J$ (on note $F$ la fibre de $J$) sur une courbe de genre $g\geq2$ avec de plus l'injection au niveau des groupes fondamentaux : $\hg{1}{F}\hookrightarrow\hg{1}{X}$. Le lemme de Castelnuevo-de Franchis (voir par exemple \cite{Bo78}) montre que $\XX_t$ admet encore une telle fibration
$$J_t:\XX_t\To C_t$$
avec $g(C_t)\geq2$ et $\hg{1}{F_t}\hookrightarrow\hg{1}{X}$ (où $F_t$, la fibre de $J_t$, est encore un tore). Le lemme \ref{additivité} montre alors que
$$\gd{\XX_t}\geq \gd{F_t}+\gd{C_t}=3$$
et on a encore :
$$\Forall{t}{\DD}\gd{\XX_t}=3=\gd{X}.$$
\vspace*{0.3cm}

\noindent\textbf{cas 3 : \mathversion{bold}$\gd{X}=3$ et $\kod{X}=2$}\mathversion{normal}\\
Là encore on va montrer que la structure de $X$ (en fait celle de son application d'Albanese) donnée par la proposition \ref{gd3/kappa2} se déforme. Considérons pour cela l'application d'Albanese relative de la famille :
$$ \xymatrix{
\XX \ar[rr]^{\alpha} \ar[rd]_{\pi} && \mathrm{Alb}(\XX/\DD) \ar[ld]\\
 & \DD}$$
\begin{enumerate}
\item si $\alpha_X$ est isotriviale sur une courbe elliptique avec pour fibre une surface de type $\pi_1$-général, il est est de même pour $\XX_t$ et on a encore $\gd{\XX_t}=3$.
\item si $\alpha_X$ est une fibration sur une surface de type $\pi_1$-général avec pour fibre des courbes $F$ vérifiant $\abs{\hg{1}{F}_X}=+\infty$, c'est encore le cas pour $\alpha_{\XX_t}$. Le lemme \ref{additivité-courbe} permet de conlure que $\XX_t$ est encore de type $\pi_1$-général.
\item si $\alpha d(X)=3$, on combine les propositions \ref{inv-albdim} et \ref{gdim-albdim} pour obtenir
$$\Forall{t}{\DD}3\geq \gd{\XX_t}\geq \alpha d(\XX_t)=\alpha d(X)=3.$$
\end{enumerate}
\vspace*{0.3cm}

\noindent\textbf{cas 4 : \mathversion{bold}$\gd{X}=2$ et $\kod{X}\leq 2$}\mathversion{normal}\\
Supposons que la $\Gamma$-dimension en soit pas localement constante ; la semi-continuité supérieure (générique) de la dimension de Kodaira et le théorème \ref{gd-générique=max} montrent alors que, génériquement sur $\DD$, on a :
$$\gd{\XX_t}=3\textrm{ et }\kod{\XX_t}\leq2.$$
Mais la discussion des cas 1, 2 et 3 entraîne qu'on doit alors avoir $\gd{X}=3$. Cette contradiction montre donc bien que
$$\gd{\XX_t}=2\textrm{ pour $t$ dans }\DD.$$

\noindent La liste des cas 1 à 4 étant exhaustive, la $\Gamma$-dimension est effectivement localement constante et ceci conclut la démonstration du théorème \ref{invariance non type general}.
\end{demo}

\section{Variétés vérifiant $\gd{X}=2$}\label{section gd=2}

Dans cette dernière partie, nous voulons comprendre plus précisément les variétés $X$ (kählérienne compacte de dimension 3) vérifiant $\gd{X}=2$. Pour cela, donnons nous un modèle holomorphe
$$\gamma:X\To \Gamma(X)=S$$
de la $\Gamma$-réduction de $X$ que nous supposerons \emph{net} : cela signifie qu'il existe un morphisme birationnel $u:X\To X'$ (avec $X'$ lisse) tel que tout diviseur $\gamma$-exceptionnel est également $u$-exceptionnel (voir \cite[def. 1.2 et lem. 1.3, p.507]{Ca04}). On introduit également le diviseur des multiplicités \textbf{au sens classique} de la fibration $\gamma$ comme dans \cite[§ 1.1.4, p. 508]{Ca04} :
$$\Delta=\Delta(\gamma)=\sum_{j\in J}(1-\frac{1}{m_j})\Delta_j$$
(les multiplicités $m_j$ sont calculées au sens pgcd) ainsi que le groupe fondamental orbifolde\footnote{celui ci est défini comme le quotient du groupe $\hg{1}{S\backslash\textrm{Supp}(\Delta)}$ par le sous-groupe normal engendré par les lacets $\beta_j^{m_j}$ où $\beta_j$ est un petit lacet autour de la composante $\Delta_j$.} $\hg{1}{S/\Delta}$. Le fait d'avoir choisi un modèle net de la fibration $\gamma$ a pour conséquence :
\begin{prop}[cor. 11.9 \cite{Ca07}]\label{suite exacte gred}
La fibration (nette) $\gamma:X\To S$ induit une suite exacte au niveau des groupes fondamentaux
$$1\To\hg{1}{F}_X \To \hg{1}{X}\To \hg{1}{S/\Delta}\To 1$$
où $F$ désigne la fibre générale de $\gamma$.
\end{prop}
Le fait de considérer le groupe fondamental orbifolde pallie donc le défaut d'exactitude (éventuel) de la suite
$$1\To\hg{1}{F}_X \To \hg{1}{X}\To \hg{1}{S}\To 1.$$

\noindent Nous allons maintenant étudier cette situation suivant les valeurs de $\kod{S/\Delta}$.

\subsection{$\Gamma$-réduction non de type général}

Traitons tout d'abord le cas où la surface orbifolde $(S/\Delta)$ n'est pas de type général. C'est ici que va intervenir la conjecture \ref{conjecture surface}. Reformulons de façon plus précise.
\begin{conj}\label{conjecture surface2}
Soit $(S/\Delta)$ une surface orbifolde lisse (au sens de \cite{Ca07}) avec $S$ rationnelle.
\begin{enumerate}
\item Si $-(K_S+\Delta)$ est ample, $\hg{1}{S/\Delta}$ est fini.
\item Si $K_S+\Delta$ est numériquement trivial, $\hg{1}{S/\Delta}$ est presque abélien.
\end{enumerate}
\end{conj}
Cette conjecture intervient naturellement dans l'étude que nous avons entreprise.
\begin{theo}[conditionné par la conjecture \ref{conjecture surface2}]\label{gred non type general}
Soit $\gamma:X\To (S/\Delta)$ comme ci-dessus ; on a alors :
\begin{description} 
\item[(i)] $\kod{S/\Delta}\geq0$.
\item[(ii)] Si $(S/\Delta)$ n'est pas de type général ($\kod{S/\Delta}<2$), il existe un revêtement étale fini $X'\To X$ de $X$ tel que la $\Gamma$-réduction de $X$ soit birationnellement équivalente à son application d'Albanese.
\end{description}
\end{theo}
\begin{rem}\label{rem S pas rationnelle}
Si on sait que $\kod{S/\Delta}=1$ ou bien que $S$ n'est pas rationnelle, les conclusions du théorème précédent sont vérifiées inconditionnellement.
\end{rem}
Le théorème \ref{gred non type general} va résulter de la combinaison des propositions suivantes. Dans les énoncés ci-dessous, nous nous plaçon toujours dans la situation où $(S/\Delta)$ est la base orbifolde (au sens classique) de l'application $\gamma_X$.
\begin{prop}\label{surface elliptique}
Si $(S/\Delta)$ admet une structure elliptique
$$f:(S/\Delta)\To C$$
sur une courbe $C$ (les fibres lisses $(E/\Delta)$ de $f$ sont elliptiques au sens orbifoldes), il existe un revêtement étale fini $X'\To X$ de $X$ vérifiant $\alpha d(X')=2$. En particulier, la $\Gamma$-réduction de $X'$ est birationnellement équivalente à son application d'Albanese $\alpha_{X'}$.
\end{prop}
\begin{demo}
Remarquons tout d'abord le fait suivant : la suite exacte des groupes fondamentaux (\ref{suite exacte gred}) montre que les fibres $E/\Delta$ de $f$ vérifient nécessairement
$$\textrm{Im}\left(\hg{1}{E/\Delta}\To \hg{1}{S/\Delta}\right)=+\infty.$$
En effet, dans la cas contraire, les sous-variétés $\Sigma=\gamma^{-1}(E)$ formeraient une famille de surfaces de $X$ avec $\hg{1}{\Sigma}_X$ fini (ce qui est contraire à l'hypothèse $\gd{X}=2$). Cette propriété est de plus stable par changement de base fini $C'\To C$ (on s'en convainc facilement en écrivant les diagrammes adéquats).

Commençons par montrer la proposition \ref{surface elliptique} dans le cas classique $\Delta=\emptyset$. Les arguments combinés de \cite[th. 2.3, p. 158]{MF94} et \cite[III 18.2 et 18.3, p. 132]{BPV} montrent que la fibration $S\To C$ est isotriviale : en effet, la présence de fibres (à réduction) singulières dans la fibration $S\To C$ entraînerait un isomorphisme au niveau des groupes fondamentaux (grâce à \cite[th. 2.3, p. 158]{MF94}) :
$$\pi_1(S)\overset{\sim}{\To}\pi_1(C)$$
ce qui est exclu par la remarque ci-dessus. On a donc $\chi_{top}(S)=0$ et les fibres de $S\To C$ sont au plus elliptiques multiples lisses. Les propositions \cite[III 18.2 et 18.3, p. 132]{BPV} s'appliquent et montrent que la fibration est effectivement isotriviale. En particulier, on en déduit que $\alpha d(S')=2$ pour un revêtement étale fini $S'$ de $S$. En considérant le revêtement (étale) $X'$ correspondant pour $X$, on a alors $\alpha d(X')\geq\alpha d(S')=2$ et la proposition \ref{gdim-albdim} montre que l'application d'Albanese de $X'$ est birationnellement équivalente à la $\Gamma$-réduction de $X'$.

Si $\Delta\neq\emptyset$, $S$ est réglée sur $C$ et les fibres orbifoldes de $f$ sont de la forme :
$$(\PP^1,\Delta)=(\PP^1/(2,3,6)),\,(\PP^1/(2,4,4)),\,(\PP^1/(3,3,3))\textrm{ ou }(\PP^1/(2,2,2,2))$$
(la notation ci-dessus suggère que seules les multiplicités nous intéressent). On peut alors appliquer le Branched Covering Trick \cite[th. 18.2, p. 56]{BPV} : un changement de base fini permet de se ramener au cas où les composantes irréductibles $f$-horizontales de $\Delta$ sont des sections de $S\To C$ (sans changer l'hypothèse sur les groupes fondamentaux par la remarque ci-dessus). En faisant agir les homographies de $\PP^1$ (ce qui induit une transformation birationnelle de la surface), on se ramène au cas où les points orbifoldes sont $0,1$ et $\infty$ (les sections correspondantes sont constantes) ; dans les trois premiers cas, la fibration $f$ est isotriviale (au sens orbifolde) et on conclut comme ci-dessus.

Dans le dernier cas, on aimerait construire le revêtement double ramifiant le long de $\textrm{Supp}(\Delta)$ pour se ramener à la situation $\Delta=\emptyset$. On constate alors que, quitte à ajouter une fibre au support de $\Delta$ (ce qui ne change pas la fibre orbifolde générale de $f$), on peut supposer que $\textrm{Supp}(\Delta)$ est divisible par 2 dans $\textrm{Pic}(S)$. Notons $\Delta^+$ ce diviseur (éventuellement) modifié, $S^+\To S$ le revêtement double le long de $\Delta^+$ et $S^+\To C^+$ la factorisation de Stein de la composée $S^+\To C$. Par la discussion menée ci-dessus, la fibration $S^+\To C^+$ est isotriviale ; \emph{a fortiori}, la fibration $f:(S/\Delta)\To C$ l'est également (comme on ne se soucie que des fibres lisses, le fait de considérer $\Delta$ ou $\Delta^+$ ne change rien). Ceci signifie exactement que la quatrième section est constante et conclut la démonstration.
\end{demo}
\begin{prop}\label{surface reglee}
Soit $(S/\Delta)$ une surface orbifolde dont la surface sous-jacente est réglée sur une courbe $C$ de genre $g\geq 1$ (\emph{via} un réglage $r:S\To C$). La structure du groupe fondamental de $(S/\Delta)$ est alors en partie dictée par la dimension de Kodaira :
\begin{enumerate}
\item Si $\kod{S/\Delta}=-\infty$, $\hg{1}{S/\Delta}$ est commensurable à celui d'une courbe de genre $g'\geq 1$.
\item Si $\kod{S/\Delta}=0$, $g=1$ et $\hg{1}{S/\Delta}$ est presque abélien.
\end{enumerate}
\end{prop}
\begin{demo}
On désignera par $E$ la fibre du réglage $r$. Si $\kod{S/\Delta}=-\infty$, on montre facilement que les fibres orbifoldes $(E/\Delta)$ vérifient également $\kod{E/\Delta}=-\infty$ (en appliquant par exemple les programmes des modèles log-minimaux \cite{FMK}). Ce sont donc des courbes rationnelles orbifoldes qui ont en particulier un groupe fondamental fini. La suite exacte des groupes fondamentaux (voir \cite{Ca07})
\begin{equation}\label{suite exacte orbifolde}
\hg{1}{E/\Delta}\To \hg{1}{S/\Delta}\To\hg{1}{C/D^*}\To 1
\end{equation}
(avec $D^*=\Delta^*(\Delta,r)$ diviseur des multiplicités de l'application $r:(S/\Delta)\To C$) montre que $\hg{1}{S/\Delta}$ est commensurable au groupe fondamental d'une courbe de genre au moins 1.\\
Si $\kod{S/\Delta}=0$, l'additivité orbifolde \cite{Ca04} a pour conséquence directe les faits suivants : $\kod{E/\Delta}=0$ et $\kod{C/D}=0$. Comme $0\leq \kod{C}\leq\kod{C/D}=0$, on en déduit que $D=\emptyset$ et $C$ est nécessairement une courbe elliptique.Via la fibration $r$, $(S/\Delta)$ est elliptique sur $C$ et la démonstration de la proposition \ref{surface elliptique} montre cette fibration est isotriviale (on peut aussi appliquer le corollaire 11.21 de \cite{Ca07} en remarquant que $\hg{1}{S/\Delta}$ est alors résiduellement fini comme extension de deux groupes eux-même résiduellement fini).
\end{demo}

\begin{demo}[du théorème \ref{gred non type general}]
\emph{Positivité de $\kod{S/\Delta}$ :}\\
Si la base orbifolde $(S/\Delta)$ de la fibration $\gamma$ vérifie $\kod{S/\Delta}=-\infty$ et si $S$ n'est pas rationnelle, la proposition \ref{surface reglee} nous montre que le groupe $\hg{1}{S/\Delta}$ est commensurable à celui d'une courbe ; d'après la suite exacte (\ref{suite exacte gred}), on en déduit que $\hg{1}{X}$ vérifie lui aussi cette relation de commensurabilité. Le théorème \ref{theoreme-siu} nous amène à la conclusion $\gd{X}=1$ ce qui n'est bien évidemment pas le cas. On a donc $\kod{S/\Delta}\geq0$ (si $S$ n'est pas rationnelle).\\

\emph{Cas $\kod{S/\Delta}=0$ :}\\
Dans ce cas, $\hg{1}{S/\Delta}$ est presque abélien : si $\kod{S}=0$, cela correspond à la conclusion de \cite[cor. 6.7, p. 347]{Ca4} et le cas $\kod{S}=-\infty$ (et $S$ non rationnelle) est traité par la proposition \ref{surface reglee}. La suite exacte (\ref{suite exacte gred}) montre alors que $\hg{1}{X}$ est également presque abélien. L'affirmation concernant l'équivalence de $\gamma_{X'}$ et $\alpha_{X'}$ (pour un revêtement étale fini $X'$ de $X$) est alors classique.\\

\emph{Cas $\kod{S/\Delta}=1$ :}\\
La surface $S$ est alors munie d'une fibration elliptique (celle donnée par la fibation d'Iitaka de $K_S+\Delta$). L'application de la proposition \ref{surface elliptique} conclut alors la démonstration.
\end{demo}

\subsection{Majoration du genre des fibres pour certaines applications de type général}

Dans ce paragraphe, on considère le cas restant, c'est-à-dire $\kod{S/\Delta}=2$. Commençons par la remarque suivante ($F$ désigne toujours la fibre de $\gamma$) :
\begin{prop}\label{gred de type general}
Dans la situation $\gamma:X\To S$ avec $\kod{S/\Delta}=2$, on a nécessairement :
$$\kod{X}=\kod{F}+\kod{S/\Delta}=\kod{F}+2.$$
En particulier,
\begin{enumerate}
\item $\kod{X}=-\infty$ si et seulement si $F=\PP^1$ (et dans ce cas, la fibration $\gamma$ coincide avec le quotient rationnel de $X$).
\item Si $\kod{X}\geq 0$, alors $\kod{X}\geq 2$ avec égalité si et seulement si $F$ est une courbe elliptique (et $\gamma$ est équivalente à la fibration d'Iitaka-Moishezon de $X$).
\end{enumerate}
\end{prop}
\begin{demo}
L'égalité annoncée n'est rien d'autre que l'additivité orbifolde \cite[cor. 4.7, p. 566]{Ca04} et les conséquences 1 et 2 sont claires.
\end{demo}

On constate donc que si $X$ n'est pas de type général, les fibres de l'application $\gamma$ sont des courbes d'un type particulier. En revanche, lorsque $\kod{X}=3$, il n'y a plus aucune restriction sur le type de courbes apparaissant comme fibre de $\gamma$. Nous allons cependant montrer que, dans cette situation, leur genre est contrôlé par le volume canonique de $X$. Rappelons que le volume d'une variété $Y$ de dimension $n$ est défini par :
$$\vol{Y}=\limsup_{m\to+\infty}\frac{h^0(Y,mK_Y)}{m^n/n!}.$$
\begin{theo}\label{majoration volume}
Soit $n\geq3$ un entier. Il existe une constante $C_n>0$ vérifiant la propriété suivante : pour toute fibration $f:Y\To Z$ de type général vers une \textbf{surface} $Z$ avec $\dimm{Y}=n$ et $Y$ de type général, le volume de la fibre générale $F$ de $f$ satisfait à l'inégalité :
$$\vol{F}\leq C_n\vol{Y}.$$
\end{theo}
\noindent Une courbe de genre $g$ a pour volume $2g-2$ et on en déduit donc :
\begin{cor}\label{majoration genre}
Dans la situation $\gamma:X\To S$ avec $\kod{S/\Delta}=2$ et $\kod{X}=3$, il existe une constante $C>0$ (indépendante de $X$) telle que :
$$g(F)\leq C\vol{X}.$$
\end{cor}

La démonstration du théorème \ref{majoration volume} repose sur une inégalité vérifiée par le volume dans les fibrations de type général (inégalité généralisant le résultat de \cite{Ka07}).
\begin{lem}\label{volume-kawamata-orbifolde}
Soit $f:Y\To Z$ une fibration de fibre générale $F$ et $\Delta$ le diviseur (classique ou non-classique) des multiplicités de $f$. Si $Y$ et $f$ sont de type général, on a alors :
$$\frac{\vol{Y}}{n!}\geq \frac{\vol{Z/\Delta}}{d!}\frac{\vol{F}}{(n-d)!}$$
avec $n=\dimm{Y}$ et $d=\dimm{Z}$.
\end{lem}
\begin{demo}
On fixe $H$ et $L$ des diviseurs amples sur $Z$ et $F$ ainsi que des entiers $m_1>m_0$ tels que :
\begin{itemize}
\item $m_0K_{(Z/\Delta)}$ et $m_1K_{(Z/\Delta)}$ soient entiers et $m_0K_{(Z/\Delta)}-H$ effectif,
\item $m_1K_F-L\geq 0$ et $\vol{L/m_1}>\vol{F}-\epsilon$.
\end{itemize}
On peut réaliser un tel choix quitte à modifier $Z$ grâce au théorème d'approximation de Fujita \cite{Fu94}. On décompose alors le fibré $kmm_1K_Y$ (pour des valeurs de $m$ et $k$ que l'on fixera dans la suite) de la façon suivante :
\begin{align*}
kmm_1K_Y&=kmm_1K_{Y\vert(Z/\Delta)}+km(m_1-m_0)f^*K_{(Z/\Delta)}\\
&\quad+kmf^*(m_0K_{(Z/\Delta)}-H)+kmf^*H.
\end{align*}
En utilisant la formule de projection et le fait que $m_0K_{(Z/\Delta)}-H$ est effectif, on obtient l'inégalité
\begin{equation}\label{fp1}
p_{kmm_1}(Y)\geq p_{k(m_1-m_0)m}(Z) h^0(Z,S^mS^k(f_*(m_1K_{Y\vert(Z/\Delta)})\otimes\mathcal{O}_Z(H))).
\end{equation}
Or, d'après \cite[th. 4.11, p. 567]{Ca04}, le faisceau $\mathcal{F}=f_*(m_1K_{Y\vert(Z/\Delta)})$ est faiblement positif et il existe donc un entier $k$ tel que $S^k(\mathcal{F}\otimes\mathcal{O}_Z(H))$ soit engendré par ses sections globales au point générique de $Y$. La dimension de l'espace des sections globales est donc au moins égale au rang de la fibre du faisceau en un point générique de $Y$ et on a donc pour $m$ assez grand :
\begin{align}\label{fp2}
h^0(Z,S^mS^k(\mathcal{F}\otimes\mathcal{O}_Z(H)))&\geq
rang\left( S^mS^k(\mathcal{F}\otimes\mathcal{O}_Z(H))\right)\nonumber\\
&\geq \dimm{S^mS^k H^0(F,m_1K_F)}\nonumber\\
&\geq\dimm{S^mS^k H^0(F,L)}\nonumber\\
&\geq h^0(F,kmL).
\end{align}
En reportant (\ref{fp2}) dans (\ref{fp1}) et en utilisant le fait que le volume de $L/m_1$ approche celui de $F$, on en déduit :
$$h^0(Y,kmm_1K_Y)\geq (\vol{Z/\Delta}-\epsilon)\frac{(k(m_1-m_0)m)^d}{d!}(\vol{F}-2\epsilon)\frac{(km_1m)^{n-d}}{(n-d)!}.$$
En choisissant $m_1$ assez grand, cette inégalité devient :
$$h^0(Y,kmm_1K_Y)\geq (\vol{Z/\Delta}-2\epsilon)(\vol{F}-2\epsilon)\frac{(km_1m)^n}{d!(n-d)!}.$$
qui, en passant à la limite, donne la minoration annoncée sur le volume de $Y$.
\end{demo}

\begin{demo}[du théorème \ref{majoration volume}]
Soit donc $f:Y\To Z$ une fibration de type général (avec $Y$ de type général) vers une surface $Z$. D'après le lemme précédent, nous disposons de l'inégalité :
$$\vol{Y}\geq \vol{F}\vol{Z/\Delta}\frac{n!}{2(n-2)!}.$$
On peut appliquer \cite{AM04} : il existe une constante universelle $\epsilon_2>0$ qui minore le volume de toute surface orbifolde (\emph{klt}) de type général :
$$\vol{Z/\Delta}\geq\epsilon_2>0.$$
L'inégalité ci-dessus se réécrit alors :
$$\vol{F}\leq C_n\vol{Y}$$
avec $C_n=\frac{2}{n(n-1)\epsilon_2}$.
\end{demo}
\begin{rem}
La minoration universelle du volume pour les orbifoldes de type général n'est malheureusement connue que dans le cas des surfaces (et comme le lemme \ref{volume-kawamata-orbifolde} est valable en toute dimension, c'est l'unique raison pour laquelle on ne considère que des fibrations vers des surfaces dans l'énoncé du théorème \ref{majoration volume}).
\end{rem}

\subsection{Ouverture de la condition $\gd{X}=3$}

Pour finir, nous donnons la démonstration du théorème \ref{ouverture gd=3}. Soit donc $\pi:\XX\To \DD$ une famille de variétés kählériennes de dimension 3 dont la fibre centrale $X=\XX_0$ vérifie $\gd{X}=3$. Au vu du théorème \ref{invariance non type general}, on peut de plus supposer que $X$ est de type général. Nous allons tout d'abord nous ramener au cas d'une famille projective grâce aux résultats de \cite[section 12.5, p. 659-663]{KM93} sur l'invariance des plurigenres en dimension 3 :
\begin{theo}\label{type-général}
Soit $\pi:\XX\To B$ une famille de variétés kählériennes de dimension 3. Si $\kod{\XX_0}=3$ pour au moins un point $0\in B$, on a
$$\Forall{b}{B}\kod{\XX_b}=3$$
et la famille $\pi:\XX\To B$ est biméromorphe (au dessus de $B$) à une famille projective.
\end{theo}
Nous pouvons donc (et c'est ce que nous ferons) supposer que la famille $\pi:\XX\To \DD$ est projective et de type général.\\

On sait alors (théorème \ref{gd-générique=max}) que, pour $t$ général dans $\DD$, $X_t$ vérifie $\gd{X_t}=3$. Si la $\Gamma$-dimension n'est pas constante au voisinage de $0\in\DD$, il existe $(t_k)_{k\geq 1}$ une suite de points de $\DD$ convergeant vers 0 telle que
$\forall k\geq1,\,\gd{X_k}=2$
(pour alléger les notations, on a posé : $X_k=X_{t_k}$). L'étude menée dans les paragraphes ci-dessus nous permet d'obtenir le résultat suivant (et donc le théorème \ref{ouverture gd=3}).
\begin{theo}\label{ouverture gd3-kappa3}
Soit $\pi:\XX\To \DD$ une famille projective de variétés de type général de dimension 3 dont la fibre centrale vérifie $\gd{X}=3$.
\begin{enumerate}
\item si il existe une suite $(t_k)_{k\geq 1}$ comme ci-dessus, la base orbifolde de la $\Gamma$-réduction de $X_k$ vérifie alors (pour $k$ assez grand) :
$$\kod{\Gamma(X_k)/\Delta^*(\gamma_k)}\leq 0.$$
\item \textbf{sous réserve de la validité de la conjecture \ref{conjecture surface2}}, la situation ci-dessus ne peut se produire et la $\Gamma$-dimension est constante au voisinage de 0 :
$$\gd{\XX_t}=3\textrm{ pour }t\textrm{ dans un voisinage de }0.$$
\end{enumerate}
\end{theo}
\begin{demo}
Remarquons tout d'abord que la possibilité :
$$\kod{\Gamma(X_k)/\Delta^*(\gamma_k)}=1$$
est exclue. En effet, d'après le théorème \ref{gred non type general} et la remarque \ref{rem S pas rationnelle}, cette hypothèse entraine que la $\Gamma$-réduction de $X_k$ est donnée (après revêtement étale fini et transformation birationnelle adéquats) par l'application d'Albanese $\alpha_{X_k}$. Or, en considérant l'application d'Albanese de la famille $\XX$, on en déduirait que la fibre centrale vérifierait $\gd{X}=2$.\\

Supposons donc que les fibres $X_k$ aient une $\Gamma$-réduction de type général :
$$\kod{\Gamma(X_k)/\Delta^*(\gamma_k)}=2.$$
Le théorème \ref{majoration volume} (ou son corollaire \ref{majoration genre}) s'applique et nous assure l'existence d'une constante universelle (numérique) $C>0$ telle que :
$$\forall k\geq1,\,g_k\leq C\vol{X_k},$$
où $g_k$ désigne le genre de la fibre générale de $\gamma_k$. Or, d'après \cite{S98}, le volume est constant dans la famille (supposée projective) et le genre des fibres générales des applications $\gamma_k$ est donc majoré par une constante indépendante de $k$. Comme le volume des fibres de $\gamma_k$ est borné et comme celles-ci restent dans un compact $\pi^{-1}(K)$ (où $K$ est par exemple un disque fermé contenant les points $t_k$), ces courbes sont donc contenues dans un compact $\mathcal{K}$ de $\mathcal{C}_1(\XX/\DD)$ (\cite{Li78}) ; en particulier, $\mathcal{K}$ ne rencontre qu'un nombre fini de composantes irréductibles de $\mathcal{C}_1(\XX/\DD)$. Quitte à extraire une sous suite de $(t_k)_k$, on peut donc trouver une composante irréductible $\mathfrak{C}$ de $\mathcal{C}_1(\XX/\DD)$ qui contient les fibres des applications $\gamma_k$. On peut alors reprendre les arguments donnés au cours de la démonstration du théorème \ref{gd-générique=max}. En effet, en passant à la limite, on peut déformer les fibres des applications $\gamma_k$ en une famille de courbes de $X$. Une composante irréductible de cette famille doit encore être couvrante sur $X$ et on conclut comme dans la démonstration du théorème \ref{gd-générique=max} : les courbes de cette famille doivent être contenues dans les fibres de la $\Gamma$-réduction de $X$ et ainsi $\gd{X}\leq 2$. Comme on a supposé $\gd{X}=3$, on en déduit bien que :
$$\textrm{pour $k$ assez grand, }\kod{\Gamma(X_k)/\Delta^*(\gamma_k)}\leq 0.$$

Si maintenant on suppose démontrée la conjecture \ref{conjecture surface2}, on sait (théorème \ref{gred non type general}) que, sous l'hypothèse
$$\kod{\Gamma(X_k)/\Delta^*(\gamma_k)}\leq 0,$$
la $\Gamma$-réduction des $X_k$ est donnée (après revêtement étale fini et transformation birationnelle) par l'application d'Albanese. A nouveau, ceci entrainerait en particulier $\gd{X}=2$ et on en déduit donc qu'une telle suite de points $(t_k)_{k\geq 1}$ ne peut exister :
$$\gd{X_t}=3\textrm{ pour tout $t$ dans un voisinage de }0.$$
\end{demo}

\section*{Remerciements}

Je tiens ici à exprimer ma plus profonde gratitude envers Frédéric Campana. Les nombreuses discussions que j'ai eue avec lui m'ont permis d'appréhender progressivement la $\Gamma$-réduction et ont largement influencé le présent travail. 

\newpage

%Bibliographie---------------------------------------------------------------
\bibliographystyle{amsalpha}
\bibliography{myref}

%Signature---------------------------------------------------------------------
\vspace*{0.5cm}
\begin{flushright}
\begin{minipage}{5cm}
Beno\^it CLAUDON\\
Universit\'e Nancy 1\\
Institut Elie Cartan\\
BP 239\\
54 506 Vandoeuvre-l\`es-Nancy\\
Cedex (France)

\vspace*{0.3cm}
\noindent Benoit.Claudon@iecn.u-nancy.fr
\end{minipage}
\end{flushright}

\end{document}